\documentclass{amsart}

\usepackage{xypic}
\input xy
\xyoption{all}
\usepackage{epsfig}
\usepackage{amsthm}
\usepackage{amssymb}
\usepackage{amsmath}
\usepackage{amscd}
\usepackage{color}
\usepackage[T1]{fontenc}

\newcommand{\bg}{\begin{equation}}
\newcommand{\ed}{\end{equation}}
\newcommand{\bga}{\begin{eqnarray}}
\newcommand{\eda}{\end{eqnarray}}
\newcommand{\pf}{\textbf{Proof:\ }}

\def\cbdu{\par{\raggedleft$\Box$\par}}

\newtheorem {Theorem}  {Theorem}

\numberwithin{Theorem}{section}

\newtheorem {Lemma}[Theorem]  {Lemma}

\theoremstyle{definition}
\newtheorem{Definition}[Theorem]{Definition}
\theoremstyle{remark}
\newtheorem{Remark}[Theorem]{\bf Remark}

\expandafter\chardef\csname pre amssym.def
at\endcsname=\the\catcode`\@ \catcode`\@=11
\def\undefine#1{\let#1\undefined}
\def\newsymbol#1#2#3#4#5{\let\next@\relax
 \ifnum#2=\@ne\let\next@\msafam@\else
 \ifnum#2=\tw@\let\next@\msbfam@\fi\fi
 \mathchardef#1="#3\next@#4#5}
\def\mathhexbox@#1#2#3{\relax
 \ifmmode\mathpalette{}{\m@th\mathchar"#1#2#3}%
 \else\leavevmode\hbox{$\m@th\mathchar"#1#2#3$}\fi}
\def\hexnumber@#1{\ifcase#1 0\or 1\or 2\or 3\or 4\or 5\or 6\or 7\or 8\or
 9\or A\or B\or C\or D\or E\or F\fi}

\font\teneufm=eufm10 \font\seveneufm=eufm7 \font\fiveeufm=eufm5
\newfam\eufmfam
\textfont\eufmfam=\teneufm \scriptfont\eufmfam=\seveneufm
\scriptscriptfont\eufmfam=\fiveeufm

\catcode`\@=\csname pre amssym.def at\endcsname

\newcounter{remark}
\newcommand{\R}{\mathbf{R}}

\renewcommand{\div}{\mbox{div}}

\def \ls {{\lambda_q^{2s}}}

\def  \R   {{\mathbb R}}

\def  \T   {{\mathbb T}}

\def  \12  {{\frac{1}{2}}}

\def\build#1_#2^#3{\mathrel{\mathop{\kern 0pt#1}\limits_{#2}^{#3}}}

\begin{document}

\title[Regularity criteria for Hall-MHD]{Regularity criterion for the 3D Hall-magneto-hydrodynamics}


\author [Mimi Dai]{Mimi Dai}
\address{Department of Mathematics, Stat. and Comp.Sci., University of Illinois Chicago, Chicago, IL 60607,USA}
\email{mdai@uic.edu}





\begin{abstract}
This paper studies the regularity problem for the 3D incompressible resistive viscous Hall-magneto-hydrodynamic (Hall-MHD) system. The Kolmogorov 41 phenomenological theory of turbulence \cite{K41} predicts that there exists a critical wavenumber above which the high frequency part is dominated by the dissipation term in the fluid equation.  Inspired by this idea, we apply an
approach of splitting the wavenumber combined with an estimate of the energy flux to obtain a new regularity criterion. The regularity condition presented here is weaker than conditions in the existing criteria (Prodi-Serrin type criteria) for the 3D Hall-MHD system.

\bigskip


\hspace{0.02cm}CLASSIFICATION CODE: 76D03, 35Q35.
\end{abstract}

\maketitle

\section{Introduction}

In this paper we consider the three dimensional incompressible resistive viscous
Hall-magneto-hydrodynamics (Hall-MHD) system:
\begin{equation}\label{HMHD}
\begin{split}
u_t+u\cdot\nabla u-b\cdot\nabla b+\nabla p=\nu\Delta u,\\
b_t+u\cdot\nabla b-b\cdot\nabla u+\nabla\times((\nabla\times b)\times b)=\mu\Delta b,\\
\nabla \cdot u=0, 
\end{split}
\end{equation}
with the initial conditions
\begin{equation}
u(x,0)=u_0(x),\qquad b(x,0)=b_0(x), \qquad \nabla\cdot u_0=\nabla\cdot b_0=0,
\end{equation}
where $x\in\mathbb{R}^3$, $t\geq 0$, $u$ is the fluid velocity, $p$ is the fluid pressure and $b$
is the magnetic field. The parameter $\nu$ denotes the kinematic viscosity coefficient of the fluid and $\mu$ denotes the reciprocal of the magnetic Reynolds number. In this paper, we assume $\nu>0$ and $\mu>0$. Note that the divergence free condition for the magnetic field $b$ is propagated by the second equation in (\ref{HMHD}) if $\nabla\cdot b_0=0$, see \cite{CL}.  One obvious difference with the usual MHD system is that the Hall-MHD system has the Hall term $\nabla\times((\nabla\times b)\times b)$ due to the happening of the magnetic reconnection when the magnetic shear is large. For the physical background of the magnetic reconnection and the Hall-MHD, we refer the readers to \cite{For, Li, SC} and references therein.

The Hall-MHD system was derived in a mathematically rigorous way by Acheritogaray, Degond, Frouvelle and Liu \cite{ADFL}, where the global existence of weak solutions in the periodic domain was
obtained. The global existence of weak solutions in the whole space $\R^3$ and the local well-posednes of classical solution were established by Chae, Degond, and Liu \cite{CDL}. The authors also obtained a blow-up criterion and the global existence of smooth solution for small initial data. Later, both the blow-up criterion and the the small data results were refined by Chae and Lee \cite{CL}. In particular, the authors proved that if a regular solution $(u,b)$ on $[0,T)$ satisfies
\begin{equation}\label{cl-criterion1}
u\in L^q(0,T;L^p(\R^3)) \qquad \mbox { and } \qquad \nabla b\in L^\gamma(0,T;L^\beta(\R^3))
\end{equation}
with 
\begin{equation}\label{cl-criterion2}
\frac3p+\frac2q\leq 1, \qquad \frac3\beta+\frac2\gamma\leq 1 \qquad \mbox { and } \qquad p,\beta\in (3, \infty]
\end{equation}
then the regular solution can be extended beyond time $T$. In the limit case $p=\beta=\infty$, it is also shown that if 
\begin{equation}\label{cl-bmo}
u, \nabla b\in L^2(0,T; BMO(\R^3))
\end{equation}
then the regular solution can be extended beyond time $T$, which is an improvement of the Prodi-Serrin condition \eqref{cl-criterion1}--\eqref{cl-criterion2}.

Partial regularity of weak solutions for the 3D Hall-MHD on plane was studied by Chae and Wolf \cite{CW}, who proved that the set of possible singularities of a weak solution has the space-time Hausdorff dimension at most two. Optimal temporal decay estimates for weak solutions were obtained by Chae and Schonbek \cite{CS}. Energy conservation for weak solutions of the 3D Hall-MHD system was studied by Dumas and Sueur \cite{DS}. Local well-posedness of classical solution to the Hall-MHD with fractional magnetic diffusion was obtained by Chae, Wan and Wu \cite{CWW}.


In this paper we will establish a new regularity criterion for the 3D Hall-MHD in term of a Besov norm with restriction only on low frequencies. We adapt the idea from the work of Cheskidov and Shyvdkoy \cite{CSreg} on the regularity problem for the Navier-Stokes equation and Euler's equation. This idea is originated from Kolmogorov's theory of turbulence, which predicts that there is a critical wavenumber above which the viscous term dominates. This method involves some techniques from harmonic analysis, such as the Littlewood-Paley decomposition theory, which are different from classical methods that have been widely used in this area. The method was also applied to improve regularity criteria for the Navier-Stokes equation and MHD system by Cheskidov and Dai in \cite{CD}, and for the supercritical quasi-geostrophic (SQG) equation by Dai in \cite{D}. Notice that the criteria obtained in \cite{CD, D} all improve the classical Prodi-Serrin, the BKM types of criteria and their extensions in each case. It suggests the wavenumber splitting method has certain advantage compared to the classical energy method in the study of the regularity problem. Therefore, we aim to apply the wavenumber splitting method to the Hall-MHD system and obtain weaker regularity condition.   

Remarkably, for the MHD system, a criterion only depending on the velocity was obtained in \cite{CD}. Namely, for a solution $(u,b)$ to the MHD system, if the velocity satisfies
\begin{equation}\label{criteria-mhd}
\displaystyle \limsup_{q\to\infty}\int_{\mathcal T_q}^T\|\Delta_q(\nabla\times u)\|_{L^\infty}dt<c
\end{equation}
for a small constant $c$, where $\Delta_q$ denotes Littlewood-Paley projection and $\{\mathcal T_q\}$ is a certain sequence of time with $\mathcal T_q\to T$ as $q\to\infty$, then the solution $(u,b)$ does not blow up at $t=T$.
Regarding the Hall-MHD system, due to the presence of the Hall term, it seems not possible to establish any criterion only in term of velocity. In the current paper, we will establish a criterion with conditions on both of the velocity and the magnetic field. Much effort will be devoted to estimating the energy flux from the Hall term that is the most difficult one. The main ingredient is the use of Littlewood-Paley decomposition theory and related estimates. For instance, Bony's paraproduct is used often to separate different types of interaction, and 
commutators are introduced to reveal cancellations contained in the nonlinear interaction.

Let $(u(t), b(t))$ be a weak solution of (\ref{HMHD}) on $[0,T]$.  Let $\lambda_q=2^{q}$, and $f_q=\Delta_qf$ is the Littlewood-Paley projection of $f$ (see Section \ref{sec:pre}). We define the dissipation wavenumber with respect to $u$ and $b$ as
\begin{equation}\label{lb}
\begin{split}
\Lambda_1(t)&=\min \left\{\lambda_q\geq1:\lambda_p^{-1}\|u_p(t)\|_\infty<c_0m, \forall p>q \right\},\\
\Lambda_2(t)&=\min \left\{\lambda_q\geq1:\lambda_{p-q}^\delta\|b_p(t)\|_\infty<c_0m, \forall p>q \right\},
\end{split}
\end{equation}
where $c_0$ is an absolute constant which will be determined later, $\lambda_{p-q}^\delta$ represents a kernel with $\delta\geq s>0$, and 
\[m=\min\{\nu,\mu\}.\] 
Let $Q_1(t), Q_2(t)\in\mathbb N$ be such that $\lambda_{Q_1(t)}=\Lambda_1(t)$ and $\lambda_{Q_2(t)}=\Lambda_2(t)$. It follows immediately that
\bg\label{Q}
\|u_{Q_1(t)}(t)\|_\infty\geq c_0m\Lambda_1(t),  \mbox { and } \quad \|\nabla b_{Q_2(t)}(t)\|_\infty\geq c_0m\Lambda_2(t).
\ed
provided $1<\Lambda_1(t), \Lambda_2(t)<\infty$; and 
\bg\notag
\|u_p(t)\|_\infty< c_0m \lambda_p \ \ \mbox { if } \ \ p>Q_1,  \mbox { and } \quad \|b_p(t)\|_\infty< c_0m \lambda_p^{-\delta} \Lambda_2^\delta(t)<c_0m  \ \ \mbox { if } \ \ p>Q_2.
\ed
 Define the function 
\begin{equation}\label{f}
f(t)=\| u_{\leq Q_1(t)}(t)\|_{ B^1_{\infty,\infty}}+ \Lambda_2(t) \|b_{\leq Q_2(t)}(t)\|_{ B^1_{\infty,\infty}},
\end{equation}
where $u_{\leq}$ and $b_{\leq}$ denote the functions restricted on low frequency part (see Section \ref{sec:pre}). Notice that $\|\nabla b_{\leq Q_2(t)}(t)\|_{ B^1_{\infty,\infty}}$ has the same scaling as $\Lambda_2(t) \|b_{\leq Q_2(t)}(t)\|_{ B^1_{\infty,\infty}}$ and is bounded by the later one.
Our main result states as follows.

\begin{Theorem}\label{thm}
Let $(u, b)$ be a weak solution to (\ref{HMHD}) on $[0,T]$. Assume that $(u(t), b(t))$ is regular on $[0,T)$, and $f\in L^1(0,T)$, i.e.
\bg\label{criteria}
\int_0^T\left(\|u_{\leq Q_1(t)}(t)\|_{ B^1_{\infty,\infty}}+ \Lambda_2(t)\|b_{\leq Q_2(t)}(t)\|_{ B^1_{\infty,\infty}}\right)\, dt<\infty.
\ed
Then $(u(t), b(t))$ is regular on $[0,T]$.
\end{Theorem}

\begin{Remark}
It will be shown in Section \ref{sec:reg} that condition \eqref{criteria} is weaker than \eqref{cl-criterion1}--\eqref{cl-criterion2} and \eqref{cl-bmo}. 
\end{Remark}

One may expect to establish a criterion analogous to (\ref{criteria-mhd}) in term of velocity and magnetic field. Concerning the length of the paper, we leave the detail for the readers who may be interested in it.

The rest of the paper is organized as follows: in Section \ref{sec:pre} we introduce some notations, recall the Littlewood-Paley decomposition theory briefly, and establish some auxiliary estimates to handle the Hall term; Section \ref{sec:reg} is devoted to proving Theorem \ref{thm}.

\bigskip

\section{Preliminaries}
\label{sec:pre}

\subsection{Notation}
\label{sec:notation}
We denote by $A\lesssim B$ an estimate of the form $A\leq C B$ with
some absolute constant $C$, and by $A\sim B$ an estimate of the form $C_1
B\leq A\leq C_2 B$ with some absolute constants $C_1$, $C_2$. We also write
 $\|\cdot\|_p=\|\cdot\|_{L^p}$, and $(\cdot, \cdot)$ stands for the $L^2$-inner product.

\subsection{Littlewood-Paley decomposition}
\label{sec:LPD}
The techniques presented in this paper rely strongly on the Littlewood-Paley decomposition. Thus we recall the Littlewood-Paley decomposition theory briefly. For a more detailed description on this theory we refer the readers to the books by Bahouri, Chemin and Danchin \cite{BCD} and Grafakos \cite{Gr}. 

Let $\mathcal F$ and $\mathcal F^{-1}$ denote the Fourier transform and inverse Fourier transform, respectively. Define $\lambda_q=2^q$ for integers $q$. A nonnegative radial function $\chi\in C_0^\infty(\R^n)$ is chosen such that 
\begin{equation}\notag
\chi(\xi)=
\begin{cases}
1, \ \ \mbox { for } |\xi|\leq\frac{3}{4}\\
0, \ \ \mbox { for } |\xi|\geq 1.
\end{cases}
\end{equation}
Let 
\bg\notag
\varphi(\xi)=\chi(\frac{\xi}{2})-\chi(\xi)
\ed
and
\begin{equation}\notag
\varphi_q(\xi)=
\begin{cases}
\varphi(\lambda_q^{-1}\xi)  \ \ \ \mbox { for } q\geq 0,\\
\chi(\xi) \ \ \ \mbox { for } q=-1.
\end{cases}
\end{equation}
For a tempered distribution vector field $u$ we define the Littlewood-Paley projection
\begin{equation}\notag
\begin{split}
&h=\mathcal F^{-1}\varphi, \qquad \tilde h=\mathcal F^{-1}\chi,\\
&u_q:=\Delta_qu=\mathcal F^{-1}(\varphi(\lambda_q^{-1}\xi)\mathcal Fu)=\lambda_q^n\int h(\lambda_qy)u(x-y)dy,  \qquad \mbox { for }  q\geq 0,\\
& u_{-1}=\mathcal F^{-1}(\chi(\xi)\mathcal Fu)=\int \tilde h(y)u(x-y)dy.
\end{split}
\end{equation}
By the Littlewood-Paley theory, the following identity
\bg\notag
u=\sum_{q=-1}^\infty u_q
\ed
holds in the distribution sense. Essentially the sequence of the smooth functions $\varphi_q$ forms a dyadic partition of the unit. To simplify the notation, we denote
\bg\notag
u_{\leq Q}=\sum_{q=-1}^Qu_q, \qquad u_{(Q, N]}=\sum_{p=Q+1}^N u_p, \qquad \tilde u_q=\sum_{|p-q|\leq 1}u_p.
\ed

\begin{Definition}
A tempered distribution $u$ belongs to the Besov space $ B_{p, \infty}^{s}$ if and only if
$$
\|u\|_{ B_{p, \infty}^{s}}=\sup_{q\geq-1}\lambda_q^s\|u_q\|_p<\infty.
$$
\end{Definition}
We also note that, 
\[
  \|u\|_{\dot H^s} \sim \left(\sum_{q=-1}^\infty\lambda_q^{2s}\|u_q\|_2^2\right)^{1/2}
\]
for each $u \in  \dot H^s$ and $s\in\R$.

We recall Bernstein's inequality for the dyadic blocks of the Littlewood-Paley decomposition in the following.
\begin{Lemma}\label{le:bern} (See \cite{L}.)
Let $n$ be the space dimension and $r\geq s\geq 1$. Then for all tempered distributions $u$, 
\bg\label{Bern}
\|u_q\|_{r}\lesssim \lambda_q^{n(\frac{1}{s}-\frac{1}{r})}\|u_q\|_{s}.
\ed
\end{Lemma}

\bigskip

\subsection{Definition of solutions}
\label{sec:sol}
We recall some classical definitions of weak and regular solutions.
\begin{Definition}\label{def:weak}
A weak solution of (\ref{HMHD}) on $[0,T]$ (or $[0, \infty)$ if $T=\infty$) is a pair of functions $(u, b)$ in the class 
$$
u, b \in C_w([0,T]; L^2(\mathbb R^3))\cap L^2(0,T; H^1(\mathbb R^3)), 
$$
with $u(0)=u_0, b(0)=b_0$, satisfying (\ref{HMHD}) in the distribution sense;
moreover, the following energy inequality 
\begin{equation}\notag
\begin{split}
&\|u(t)\|_2^2+\|b(t)\|_2^2+2\nu\int_{t_0}^t\|\nabla u(s)\|_2^2ds
+2\mu\int_{t_0}^t\|\nabla b(s)\|_2^2ds\\
\leq &\|u(t_0)\|_2^2+\|b(t_0)\|_2^2
\end{split}
\end{equation}
is satisfied for almost all $t_0\in(0, T)$ and all $t\in(t_0, T]$.
\end{Definition}

\begin{Lemma}\label{le-reg} (See \cite{CL}.)
A weak solution $(u, b)$ of (\ref{HMHD}) is regular on a time interval $\mathcal I$ if $\|u(t)\|_{H^s}$ and $\|b(t)\|_{H^s}$ are continuous on $\mathcal I$ for some $s> \frac{5}{2}$.
\end{Lemma}

\subsection{Bony's paraproduct and commutator}
\label{sec-para}

Bony's paraproduct formula 
\begin{equation}\label{Bony}
\begin{split}
\Delta_q(u\cdot\nabla v)=&\sum_{|q-p|\leq 2}\Delta_q(u_{\leq{p-2}}\cdot\nabla v_p)+
\sum_{|q-p|\leq 2}\Delta_q(u_{p}\cdot\nabla v_{\leq{p-2}})\\
&+\sum_{p\geq q-2} \Delta_q(\tilde u_p \cdot\nabla v_p),
\end{split}
\end{equation}
will be used constantly to decompose the nonlinear terms in energy estimate. 
We will also use the notation of the commutator
\begin{equation} \label{commudef}
[\Delta_q, u_{\leq{p-2}}\cdot\nabla]v_p:=\Delta_q(u_{\leq{p-2}}\cdot\nabla v_p)-u_{\leq{p-2}}\cdot\nabla \Delta_qv_p.
\end{equation}
\begin{Lemma}\label{le-commu}
The commutator satisfies the following estimate, for any $1<r<\infty$
\[\|[\Delta_q,u_{\leq{p-2}}\cdot\nabla] v_p\|_{r}\lesssim \|\nabla u_{\leq p-2}\|_\infty\|v_p\|_{r}.\]
\end{Lemma}
\pf
It follows from the definition of $\Delta_q$ that
\begin{equation}\notag
\begin{split}
[\Delta_q,u_{\leq{p-2}}\cdot\nabla] v_p=&\int_{\T^2}\lambda_q^3h(\lambda_q(x-y))\left(u_{\leq p-2}(x)-u_{\leq p-2}(y)\right)\nabla_y v_p(y)\,dy\\
=&-\int_{\T^2}\lambda_q^3\nabla_y h(\lambda_q(x-y))\left(u_{\leq p-2}(x)-u_{\leq p-2}(y)\right) v_p(y)\,dy\\
=&-\int_{\T^2}\lambda_q^3|x-y|\nabla_y h(\lambda_q(x-y))\frac{u_{\leq p-2}(x)-u_{\leq p-2}(y)}{|x-y|} v_p(y)\,dy
\end{split}
\end{equation}
where we used the integration by parts and the fact $\div\, u_{\leq p-2}=0$.
Thus, by Young's inequality, for any $r>1$,
\begin{equation}\notag
\begin{split}
&\|[\Delta_q,u_{\leq{p-2}}\cdot\nabla] v_p\|_{r}\\
&\lesssim \|\nabla u_{\leq p-2}\|_\infty\|v_p\|_{r}\left|\int_{\T^2}|z||\nabla h(z)|\, dz\right|\\
&\lesssim \|\nabla u_{\leq p-2}\|_\infty\|v_p\|_{r}.
\end{split}
\end{equation}
\cbdu

\subsection{Auxiliary estimates}
\label{sec-commu}

To handle the Hall term $\nabla \times((\nabla \times b)\times b)$, we introduce a commutator and some other estimates that follow from some elementary vector calculus.
Define the commutator for vector valued functions $F$ and $G$,
\begin{equation}\label{comm-v1}
[\Delta_q,F\times\nabla\times]G=\Delta_q(F\times(\nabla\times G))-F\times(\nabla\times G_q);
\end{equation}
\begin{equation}\label{comm-v2}
[\Delta_q, (\nabla\times F)\times]G=\Delta_q((\nabla\times F)\times G)-(\nabla \times F)\times G_q.
\end{equation}
The commutator will be used to reveal certain cancellation from the Hall term. It satisfies the following estimates. 
\begin{Lemma}\label{le-Hall1}
Let $F$ and $G$ be vector valued functions. Assume $\nabla\cdot F=0$ and $F$, $G$ vanish at large $|x|\in \R^3$. For any $1<r<\infty$, we have
\[\|[\Delta_q,F\times\nabla\times]G\|_r\lesssim  \|\nabla F\|_\infty\|G\|_r;\]
\[\|[\Delta_q, (\nabla\times F)\times]G\|_r\lesssim  \|\nabla F\|_\infty\|G\|_r;\]
\end{Lemma}
\pf
Due to the fact $\nabla\cdot F=0$, using integration by parts yields that for the scalar function $h$, 
\begin{equation}\label{IBP}
\begin{split}
\int_{\mathbb R^3} h(\nabla \times G)\times F\,dx=&-\int_{\mathbb R^3} (\nabla h \times G)\times F\,dx+\int_{\mathbb R^3} h G \times (\nabla \times F)\,dx\\
&+\int_{\mathbb R^3} h G \cdot \nabla  F\,dx.
\end{split}
\end{equation} 
It follows from  the definition of $\Delta_q$ and (\ref{IBP}) that
\begin{equation}\notag
\begin{split}
[\Delta_q,F\times\nabla\times]G=& \Delta_q(F\times(\nabla\times G))-F\times(\nabla\times G_q)\\
=&\int_{\R^3}\lambda_q^3h(\lambda_q(x-y))(F(x)-F(y))\times\nabla_y\times G(y)\, dy\\
=& \int_{\R^3}\lambda_q^3\nabla_y h(\lambda_q(x-y))\times G(y)\times(F(x)-F(y))\, dy\\
&-\int_{\R^3} \lambda_q^3h(\lambda_q(x-y)) G(y)\times(\nabla_y\times(F(x)-F(y)))\, dy\\
&-\int_{\R^3} \lambda_q^3h(\lambda_q(x-y)) G(y)\cdot\nabla_y (F(x)-F(y))\, dy,\\
\end{split}
\end{equation}
while the first integral can be rewritten as 
\[\int_{\R^3} \lambda_q^3|x-y|\nabla_y h(\lambda_q(x-y))\times G(y)\times\frac{(F(x)-F(y))}{|x-y|}\, dy.\]
Therefore, it follows from Young's inequality that
\begin{equation}\notag
\begin{split}
\|[\Delta_q,F\times\nabla\times]G\|_r
\lesssim & \|\nabla F\|_\infty\|G\|_r\int_{\R^3}\lambda_q^3|x-y||\nabla_y h(\lambda_q(x-y))|\, dy\\
&+\|\nabla F\|_\infty\|G\|_r\int_{\R^3}\lambda_q^3|h(\lambda_q(x-y))|\, dy\\
\lesssim & \|\nabla F\|_\infty\|G\|_r.
\end{split}
\end{equation}
Another inequality in the lemma can be obtained in an analogous way.

\cbdu

\begin{Lemma}\label{le-Hall2}
Assume the vector valued functions $F$, $G$ and $H$ vanish at large $|x|\in \R^3$. For any $1\leq r_1, r_2\leq \infty$ with $\frac1{r_1}+\frac1{r_2}=1$,  we have
\[
\left|\int_{\R^3}[\Delta_q, (\nabla\times F)\times]G\cdot\nabla\times H\, dx\right|
\lesssim \|\nabla^2 F\|_\infty\|G\|_{r_1}\|H\|_{r_2}.
\]
\end{Lemma}
\pf 
The definition of $\Delta_q$ along with (\ref{comm-v2}) indicates that
\begin{equation}\notag
\begin{split}
&\int_{\R^3}[\Delta_q,(\nabla\times F)\times]G\cdot\nabla\times H\, dx\\
=&\int_{\R^3}\int_{\R^3}\lambda_q^3h(\lambda_q(x-y))\left[\nabla_x\times F(x)-\nabla_y\times F(y)\right]\times G(y)\cdot\nabla_x\times H(x)\, dy\,dx\\
=&-\int_{\R^3}\int_{\R^3}\lambda_q^3\nabla_x h(\lambda_q(x-y))\left[\nabla_x\times F(x)-\nabla_y\times F(y)\right]\times G(y)\cdot H(x)\, dy\,dx\\
&-\int_{\R^3}\int_{\R^3}\lambda_q^3h(\lambda_q(x-y))\nabla_x\times\left[\nabla_x\times F(x)-\nabla_y\times F(y)\right]\times G(y)\cdot H(x)\, dy\,dx\\
\equiv &J_1+J_2.
\end{split}
\end{equation}
It then follows  from Young's convolution inequality that, for $\frac{1}{r_1}+\frac{1}{r_2}=1$ with $1\leq r_1,r_2\leq \infty$,
\begin{equation}\notag
\begin{split}
|J_{1}|=&\left|\int_{\R^3}\int_{\R^3}\lambda_q^3|x-y|\nabla_x h(\lambda_q(x-y))\frac{\left[\nabla_x\times F(x)-\nabla_y\times F(y)\right]}{|x-y|}\times G(y)\cdot H(x)\, dy\,dx\right|\\
\lesssim &\|G\|_{r_1}\|H\|_{r_2}\left\|\frac{\left[\nabla_x\times F(x)-\nabla_y\times F(y)\right]}{|x-y|}\right\|_\infty\int_{\R^3}\lambda_q^3|x-y|\left|\nabla_x h(\lambda_q(x-y))\right|\, dy\\
\lesssim &\|\nabla^2 F\|_\infty\|G\|_{r_1}\|H\|_{r_2}.
\end{split}
\end{equation}
Applying Young's convolution inequality to $J_2$ yields
 \begin{equation}\notag
\begin{split}
|J_2|\leq &  \int_{\R^3}\int_{\R^3}\lambda_q^3\left|h(\lambda_q(x-y))[\nabla_x\times\nabla_x\times F(x)]\times G(y)\right|\, dy |H(x)|\, dx\\
\lesssim &\|\nabla^2 F\|_\infty\|G\|_{r_1}\|H\|_{r_2}\int_{\R^3}\lambda_q^3\left| h(\lambda_q(x-y))\right|\, dy\\
\lesssim &\|\nabla^2 F\|_\infty\|G\|_{r_1}\|H\|_{r_2}.
\end{split}
\end{equation}

\cbdu

\bigskip

\section{Regularity criterion}
\label{sec:reg}

\subsection{Proof of Theorem \ref{thm}}
Thanks to Lemma \ref{le-reg}, in order to prove the weak solution $(u,b)$ is regular, it is sufficient to prove that $\|u(t)\|_{H^s}+\|b(t)\|_{H^s}$ is  bounded on $[0,T)$ for some $s>\frac{5}{2}$.  
Since $(u(t), b(t))$ is regular on $[0,T)$, multiplying the equations of (\ref{HMHD}) with $\Delta^2_qu$ and $\Delta^2_qb$ respectively yields 
\begin{equation}\notag
\begin{split}
\frac{1}{2}\frac{d}{dt}\|u_q\|_2^2\leq &-\nu\|\nabla u_q\|_2^2+\int_{\R^3}\Delta_q(u\cdot\nabla u)\cdot u_q\, dx
-\int_{\R^3}\Delta_q(b\cdot\nabla b)\cdot u_q\, dx\\
\frac{1}{2}\frac{d}{dt}\|b_q\|_2^2\leq &-\mu\|\nabla b_q\|_2^2+\int_{\R^3}\Delta_q(u\cdot\nabla b)\cdot b_qdx-\int_{\R^3}\Delta_q(b\cdot\nabla u)\cdot b_qdx\\
&-\int_{\R^3}\Delta_q((\nabla\times b)\times b)\cdot \nabla\times b_q\, dx.
\end{split}
\end{equation}
Multiplying the above two inequalities by $\lambda_q^{2s}$ and adding up for all $q\geq -1$ we obtain 
\begin{equation}\label{ineq-ubq}
\begin{split}
\frac{1}{2}\frac{d}{dt}\sum_{q\geq -1}\lambda_q^{2s}\left(\|u_q\|_2^2+\|b_q\|_2^2\right)\leq &-\sum_{q\geq -1}\lambda_q^{2s}\left(\nu\|\nabla u_q\|_2^2+\mu\|\nabla b_q\|_2^2\right)\\
&+I_1+I_2+I_3+I_4+I_5,
\end{split}
\end{equation}
with
\begin{equation}\notag
\begin{split}
I_1=&\sum_{q\geq -1}\lambda_q^{2s}\int_{\R^3}\Delta_q(u\cdot\nabla u)\cdot u_q\, dx, \qquad
I_2=-\sum_{q\geq -1}\lambda_q^{2s}\int_{\R^3}\Delta_q(b\cdot\nabla b)\cdot u_q\, dx,\\
I_3=&\sum_{q\geq -1}\lambda_q^{2s}\int_{\R^3}\Delta_q(u\cdot\nabla b)\cdot b_q\, dx,\qquad
I_4=-\sum_{q\geq -1}\lambda_q^{2s}\int_{\R^3}\Delta_q(b\cdot\nabla u)\cdot b_q\, dx,\\
I_5=&-\sum_{q\geq -1}\lambda_q^{2s}\int_{\R^3}\Delta_q((\nabla\times b)\times b)\cdot \nabla\times b_q\, dx.
\end{split}
\end{equation}
The idea is to establish a Gr\"onwall's type inequality for $\|u\|_{H^s}^2+\|b\|_{H^s}^2$. The main ingredients to estimate the terms $I_1, \ldots, I_5$ include the usage of  Bony's para-product and commutators mentioned in Section \ref{sec:pre}. Typically, commutators help us to move derivatives from high frequency to low frequency terms and also reveal cancellations in the setting of the Littlewood-Paley decomposition.

Notice that the flux terms $I_1, I_2, I_3, I_4$ have been estimated in \cite{CD} where criterion only in term of velocity for the 3D MHD was obtained. The situation in this paper is different since the regularity condition will be on both of the velocity and the magnetic field. Thus, $I_2$ and $I_4$ will be estimated in a slight different way which requires less restriction on $s$, while $I_1$ and $I_3$ can be estimated the same way as in \cite{CD} (taking $r=\infty$). In the end, we will estimate $I_5$ in detail, which is the most difficult term due to the strong nonlinearity.

Also notice that the definition of $f(t)$ in \cite{CD} is different from (\ref{f}). We omit the details of computation and conclude that, for any $s>\frac12$
 \begin{equation}\label{est-i1}
|I_1|\lesssim c_0m\sum_{q\geq -1}\lambda_q^{2s+2}(\|u_q\|_2^2+\|b_q\|_2^2)+Q_1f(t)\sum_{q\geq -1}\lambda_q^{2s}(\|u_q\|_2^2+\|b_q\|_2^2).
\end{equation}
We estimate $I_2$ and $I_4$ in the following and show that cancellation occurs in $I_2+I_4$.
Using Bony's paraproduct and the commutator notation, $I_2$ is decomposed as
\begin{equation}\notag
\begin{split}
I_2=
&-\sum_{q\geq -1}\sum_{|q-p|\leq 2}\lambda_q^{2s}\int_{\R^3}\Delta_q(b_{\leq p-2}\cdot\nabla b_p)u_q\, dx\\
&-\sum_{q\geq -1}\sum_{|q-p|\leq 2}\lambda_q^{2s}\int_{\R^3}\Delta_q(b_{p}\cdot\nabla b_{\leq{p-2}})u_q\, dx\\
&-\sum_{q\geq -1}\sum_{p\geq q-2}\lambda_q^{2s}\int_{\R^3}\Delta_q(b_p\cdot\nabla\tilde b_p)u_q\, dx\\
=&I_{21}+I_{22}+I_{23},
\end{split}
\end{equation}
with 
\begin{equation}\notag
\begin{split}
I_{21}=&-\sum_{q\geq -1}\sum_{|q-p|\leq 2}\lambda_q^{2s}\int_{\R^3}[\Delta_q, b_{\leq{p-2}}\cdot\nabla] b_pu_q\, dx\\
&-\sum_{q\geq -1}\sum_{|q-p|\leq 2}\lambda_q^{2s}\int_{\R^3}b_{\leq{q-2}}\cdot\nabla \Delta_q b_p u_q\, dx\\
&-\sum_{q\geq -1}\sum_{|q-p|\leq 2}\lambda_q^{2s}\int_{\R^3}(b_{\leq{p-2}}-b_{\leq{q-2}})\cdot\nabla\Delta_qb_p u_q\, dx\\
=&I_{211}+I_{212}+I_{213}.
\end{split}
\end{equation}
We will see that the term $I_{212}$ cancels a part from $I_4$. The other terms are estimated as follows.
Splitting the summation by the wavenumber $\Lambda_2$ yields,  
\begin{equation}\notag
\begin{split}
|I_{211}|\leq &\sum_{p\geq 1}\sum_{|q-p|\leq 2}\lambda_q^{2s}\int_{\R^3}\left|[\Delta_q, b_{\leq{p-2}}\cdot\nabla] b_pu_q\right|\, dx\\
\leq & \sum_{1\leq p\leq Q_2+2}\sum_{|q-p|\leq 2}\lambda_q^{2s}\int_{\R^3}\left|[\Delta_q, b_{\leq{p-2}}\cdot\nabla] b_pu_q\right|\, dx\\
&+\sum_{p> Q_2+2}\sum_{|q-p|\leq 2}\lambda_q^{2s}\int_{\R^3}\left|[\Delta_q, b_{\leq Q_2}\cdot\nabla] b_pu_q\right|\, dx\\
&+\sum_{p> Q_2+2}\sum_{|q-p|\leq 2}\lambda_q^{2s}\int_{\R^3}\left|[\Delta_q, b_{(Q_2,p-2]}\cdot\nabla] b_pu_q\right|\, dx\\
\equiv &A_1+A_2+A_3;
\end{split}
\end{equation}
with 
\begin{equation}\notag
\begin{split}
A_1\lesssim &\sum_{1\leq p\leq Q_2+2}\sum_{|q-p|\leq 2}\lambda_q^{2s}\|\nabla b_{\leq p-2}\|_\infty\|b_p\|_2\|u_q\|_2\\
\lesssim & Q_2f(t) \sum_{1\leq p\leq Q_2+2}\|b_p\|_2\sum_{|q-p|\leq 2}\lambda_q^{2s}\|u_q\|_2\\
\lesssim & Q_2f(t) \sum_{q\geq -1}\lambda_q^{2s}\left(\|u_q\|_2^2+\|b_q\|_2^2\right);
\end{split}
\end{equation}

\begin{equation}\notag
\begin{split}
A_2\lesssim &\sum_{ p> Q_2+2}\sum_{|q-p|\leq 2}\lambda_q^{2s}\|\nabla b_{\leq Q_2}\|_\infty\|b_p\|_2\|u_q\|_2\\
\lesssim & Q_2f(t) \sum_{ p> Q_2+2}\|b_p\|_2\sum_{|q-p|\leq 2}\lambda_q^{2s}\|u_q\|_2\\
\lesssim & Q_2f(t) \sum_{q> Q_2}\lambda_q^{2s}\left(\|u_q\|_2^2+\|b_q\|_2^2\right);
\end{split}
\end{equation}

\begin{equation}\notag
\begin{split}
A_3\lesssim &\sum_{ p> Q_2+2}\sum_{|q-p|\leq 2}\lambda_q^{2s}\|\nabla b_{(Q_2,p-2]}\|_\infty\|b_p\|_2\|u_q\|_2\\
\lesssim & \sum_{ p> Q_2+2}\|b_p\|_2\sum_{|q-p|\leq 2}\lambda_q^{2s}\|u_q\|_2\sum_{Q_2<p'\leq p-2}\lambda_{p'}\|b_{p'}\|_\infty\\
\lesssim & c_0m \sum_{ p> Q_2+2}\|b_p\|_2\sum_{|q-p|\leq 2}\lambda_q^{2s}\|u_q\|_2\sum_{Q_2<p'\leq p-2}\lambda_{p'}^{2}\\
\lesssim & c_0m \sum_{ p> Q_2+2}\lambda_p^{2s}\left(\|u_p\|_2^2+\|b_p\|_2^2\right)\sum_{Q_2<p'\leq p-2}\lambda_{p'}^{2}\\
\lesssim & c_0m \sum_{ p> Q_2+2}\lambda_p^{2s+2}\left(\|u_p\|_2^2+\|b_p\|_2^2\right)\sum_{Q_2<p'\leq p-2}\lambda_{p'-p}^{2}\\
\lesssim & c_0m \sum_{ p> Q_2+2}\lambda_p^{2s+2}\left(\|u_p\|_2^2+\|b_p\|_2^2\right).
\end{split}
\end{equation}
The term $I_{213}$ is estimated as
\begin{equation}\notag
\begin{split}
|I_{213}|\leq & \sum_{q\geq -1}\sum_{|q-p|\leq 2}\lambda_q^{2s}\int_{\R^3}|(b_{\leq{p-2}}-b_{\leq{q-2}})\cdot\nabla\Delta_qb_p u_q|\, dx\\
\leq & \sum_{q>Q_1}\sum_{|q-p|\leq 2}\lambda_q^{2s+1}\|u_q\|_\infty\|b_{\leq{p-2}}-b_{\leq{q-2}}\|_2\|b_p\|_2\\
& +\sum_{-1\leq q\leq Q_1}\sum_{|q-p|\leq 2}\lambda_q^{2s+1}\|u_q\|_\infty\|b_{\leq{p-2}}-b_{\leq{q-2}}\|_2\|b_p\|_2\\
\lesssim & c_0m \sum_{q>Q_1}\lambda_q^{2s+2}\sum_{|q-p|\leq 2}\|b_{\leq{p-2}}-b_{\leq{q-2}}\|_2\|b_p\|_2\\
& +f(t)\sum_{-1\leq q\leq Q_1}\lambda_q^{2s}\sum_{|q-p|\leq 2}\|b_{\leq{p-2}}-b_{\leq{q-2}}\|_2\|b_p\|_2\\
\lesssim & c_0m \sum_{q>Q_1}\lambda_q^{2s+2}\sum_{q-3\leq p\leq q+2}\|b_p\|_2^2 +f(t)\sum_{-1\leq q\leq Q_1}\lambda_q^{2s}\|b_q\|_2^2\\
\lesssim & c_0m \sum_{q>Q_1-3}\lambda_q^{2s+2}\|b_p\|_2^2 +f(t)\sum_{-1\leq q\leq Q_1}\lambda_q^{2s}\|b_q\|_2^2.
\end{split}
\end{equation}
Notice that $I_{22}$ has the same estimate as $I_{211}$.  While $I_{23}$ is estimated as
\begin{equation}\notag
\begin{split}
|I_{23}|\lesssim &\sum_{q\geq -1}\sum_{p\geq q-2}\lambda_q^{2s}\int_{\mathbb R^3}|\Delta_q(b_p\otimes\tilde b_p)\nabla u_q|\, dx\\
\lesssim &\sum_{q> Q_1}\lambda_q^{2s+1}\|u_q\|_\infty\sum_{p\geq q-2}\|b_p\|_2^2
+\sum_{-1\leq q\leq Q_1}\lambda_q^{2s+1}\|u_q\|_\infty\sum_{p\geq q-2}\|b_p\|_2^2\\
\lesssim &c_0m\sum_{q> Q_1}\lambda_q^{2s+2}\sum_{p\geq q-2}\|b_p\|_2^2 +f(t)\sum_{-1\leq q\leq Q_1}\lambda_q^{2s}\sum_{p\geq q-2}\|b_p\|_2^2\\
\lesssim &c_0m\sum_{p> Q_1}\lambda_p^{2s+2}\|b_p\|_2^2\sum_{Q_1< q\leq p+2}\lambda_{q-p}^{2s+2}\\
&+f(t)\sum_{-1\leq p\leq Q_1}\lambda_p^{2s}\|b_p\|_2^2\sum_{-1\leq q\leq p}\lambda_{q-p}^{2s}
+f(t)\sum_{p> Q_1}\lambda_p^{2s}\|b_p\|_2^2\sum_{-1\leq q\leq Q_1}\lambda_{q-p}^{2s}\\
\lesssim &c_0m\sum_{q> Q_1}\lambda_q^{2s+2}\|b_q\|_2^2+f(t)\sum_{q\geq -1}\lambda_q^{2s}\|b_q\|_2^2.
\end{split}
\end{equation}
Therefore, we have for $s>\frac12$
\begin{equation}\label{est-i2}
|I_2|\lesssim c_0m\sum_{q\geq -1}\lambda_q^{2s+2}(\|u_q\|_2^2+\|b_q\|_2^2)+Q_2f(t)\sum_{q\geq -1}\lambda_q^{2s}(\|u_q\|_2^2+\|b_q\|_2^2).
\end{equation}

Now we estimate $I_4$ in a similar way. By Bony's paraproduct and the commutator notation, $I_4$ can be decomposed as
\begin{equation}\notag
\begin{split}
I_4=
&-\sum_{q\geq -1}\sum_{|q-p|\leq 2}\lambda_q^{2s}\int_{\R^3}\Delta_q(b_{\leq p-2}\cdot\nabla u_p)b_q\, dx\\
&-\sum_{q\geq -1}\sum_{|q-p|\leq 2}\lambda_q^{2s}\int_{\R^3}\Delta_q(b_{p}\cdot\nabla u_{\leq{p-2}})b_q\, dx\\
&-\sum_{q\geq -1}\sum_{p\geq q-2}\lambda_q^{2s}\int_{\R^3}\Delta_q(\tilde b_p\cdot\nabla u_p)b_q\, dx\\
=&I_{41}+I_{42}+I_{43},
\end{split}
\end{equation}
with 
\begin{equation}\notag
\begin{split}
I_{41}=&-\sum_{q\geq -1}\sum_{|q-p|\leq 2}\lambda_q^{2s}\int_{\R^3}[\Delta_q,b_{\leq{p-2}}\cdot\nabla] u_pb_q\, dx\\
&-\sum_{q\geq -1}\sum_{|q-p|\leq 2}\lambda_q^{2s}\int_{\R^3}b_{\leq{q-2}}\cdot\nabla\Delta_q u_p b_q\, dx\\
&-\sum_{q\geq -1}\sum_{|q-p|\leq 2}\lambda_q^{2s}\int_{\R^3}(b_{\leq{p-2}}-b_{\leq{q-2}})\cdot\nabla\Delta_qu_p b_q\, dx\\
=&I_{411}+I_{412}+I_{413}.
\end{split}
\end{equation}
As mentioned above the term $I_{412}$ cancels $I_{212}$. Indeed, we have, using integration by parts 
\begin{equation}\notag
\begin{split}
I_{212}+I_{412}=&-\sum_{q\geq -1}\sum_{|q-p|\leq 2}\lambda_q^{2s}\int_{\R^3}b_{\leq{q-2}}\cdot\nabla\Delta_q b_p (u_q+b_q)\, dx\\
&-\sum_{q\geq -1}\sum_{|q-p|\leq 2}\lambda_q^{2s}\int_{\R^3}b_{\leq{q-2}}\cdot\nabla\Delta_q u_p (b_q+u_q)\, dx\\
=&-\sum_{q\geq -1}\lambda_q^{2s}\int_{\R^3}b_{\leq{q-2}}\cdot\nabla(u_p+b_p) (u_q+b_q)\, dx\\
=&0.
\end{split}
\end{equation}
Notice that $I_{411}$ can be estimated as $I_{211}$, $I_{42}$ can be estimated as $I_{22}$, and $I_{43}$ can be estimated as a similar part from $I_{3}$, thus
\begin{equation}\notag
|I_{411}|+|I_{42}+|I_{43}|\lesssim c_0m\sum_{q\geq -1}\lambda_q^{2s+2}(\|u_q\|_2^2+\|b_q\|_2^2)+Q_1f(t)\sum_{q\geq -1}\lambda_q^{2s}(\|u_q\|_2^2+\|b_q\|_2^2).
\end{equation}
After using integration by parts, the term $I_{413}$ can be estimated similarly as $I_{213}$, hence
\begin{equation}\notag
|I_{413}|\lesssim c_0m\sum_{q\geq -1}\lambda_q^{2s+2}\|b_q\|_2^2+f(t)\sum_{q\geq -1}\lambda_q^{2s}\|b_q\|_2^2.
\end{equation} 
Thus, we obtain 
\begin{equation}\label{est-i4}
|I_4|\lesssim c_0m\sum_{q\geq -1}\lambda_q^{2s+2}(\|u_q\|_2^2+\|b_q\|_2^2)+Q_1f(t)\sum_{q\geq -1}\lambda_q^{2s}(\|u_q\|_2^2+\|b_q\|_2^2).
\end{equation} 

In the following, we focus on the estimate for $I_5=-\lambda_q^{2s}\int_{\R^3}\Delta_q((\nabla\times b)\times b)\cdot \nabla\times b_q\, dx$, which comes from the Hall term. The Hall term involves the strongest nonlinearity in the equation and thus is the most difficult term to estimate. Specifically, the local high frequency interactions accumulate to a large and hard to control  term.  Thanks to the commutators (\ref{comm-v1}) and (\ref{comm-v2}), a decomposition is applied so that  $I_{512}$, which contains the worst interaction, actually vanishes. While the other terms left in the decomposition are estimated by the auxiliary estimates established in Subsection \ref{sec-commu}.

Applying Bony's paraproduct first,  $I_5$ is decomposed as
\begin{equation}\notag
\begin{split}
I_{5}=&\sum_{q\geq-1}\sum_{|q-p|\leq 2}\ls\int_{\mathbb R^3}\Delta_q( b_{\leq p-2}\times(\nabla\times b_p))\cdot\nabla\times b_q\, dx\\
&+\sum_{q\geq-1}\sum_{|q-p|\leq 2}\ls\int_{\mathbb R^3}\Delta_q( b_{p}\times(\nabla\times b_{\leq p-2}))\cdot\nabla\times b_q\, dx\\
&+\sum_{q\geq-1}\sum_{p\geq q-2}\ls\int_{\mathbb R^3}\Delta_q( b_{p}\times(\nabla\times \tilde b_p))\cdot\nabla\times b_q\, dx\\
=&I_{51}+I_{52}+I_{53}.
\end{split}
\end{equation}
Using the commutator notation (\ref{comm-v1}), $I_{51}$ can be further decomposed as
\begin{equation}\notag
\begin{split}
I_{51}=&\sum_{q\geq-1}\sum_{|q-p|\leq 2}\ls\int_{\mathbb R^3}[\Delta_q,b_{\leq p-2}\times\nabla\times]b_p\cdot\nabla\times b_q\, dx\\
&+\sum_{q\geq-1}\ls\int_{\mathbb R^3}b_{\leq q-2}\times(\nabla\times b_q)\cdot\nabla\times b_q\, dx\\
&+\sum_{q\geq-1}\sum_{|p-q|\leq 2}\ls\int_{\mathbb R^3}(b_{\leq p-2}-b_{\leq q-2})\times(\nabla\times (b_p)_q)\cdot\nabla\times b_q\, dx\\
=&I_{511}+I_{512}+I_{513},
\end{split}
\end{equation}
where we used the fact $\sum_{|p-q|\leq 2}\Delta_qb_p=b_q$.
It is clear that $I_{512}=0$ due to property of cross product. 
While, we have
\begin{equation}\notag
\begin{split}
|I_{511}|
\leq&\sum_{1\leq p\leq Q_2+2}\sum_{|p-q|\leq 2}\lambda_q^{2s}\int_{\mathbb R^3}|[\Delta_q,b_{\leq p-2}\times\nabla\times]b_p\cdot\nabla\times b_q|\, dx\\
&+\sum_{p> Q_2+2}\sum_{|p-q|\leq 2}\lambda_q^{2s}\int_{\mathbb R^3}|[\Delta_q,b_{\leq Q_2}\times\nabla\times]b_p\cdot\nabla\times b_q|  \, dx\\
&+\sum_{p>Q_2+2}\sum_{|p-q|\leq 2}\lambda_q^{2s}  \int_{\mathbb R^3}|[\Delta_q,b_{(Q_2,p-2]}\times\nabla\times]b_p\cdot\nabla\times b_q|  \, dx \\
\equiv &A_{4}+A_{5}+A_{6}.
\end{split}
\end{equation}
By Lemma  \ref{le-Hall1}, the definition of $f(t)$ (\ref{f}), we infer
\begin{equation}\notag
\begin{split}
A_{4}\lesssim &\sum_{1\leq p\leq Q_2+2}\sum_{|p-q|\leq 2}\lambda_q^{2s+1}\|\nabla b_{\leq p-2}\|_\infty\|b_p\|_2\|b_q\|_2\\
\lesssim & Q_2\Lambda_2 \|\nabla b_{\leq Q_2}\|_\infty \sum_{1\leq p\leq Q_2+2}\|b_p\|_2\sum_{|p-q|\leq 2}\lambda_q^{2s}\|b_q\|_2\\
\lesssim & Q_2\Lambda_2 \|\nabla b_{\leq Q_2}\|_\infty\sum_{p\geq -1}\lambda_p^{2s}\|b_p\|_2^2.
\end{split}
\end{equation}
To estimate $A_5$, we deduce by using Lemma \ref{le-Hall1} and the definition of $\Lambda_2(t)$ (\ref{lb})
\begin{equation}\notag
\begin{split}
A_{5}\lesssim &\sum_{p> Q_2+2}\sum_{|p-q|\leq 2}\lambda_q^{2s+1}\|\nabla b_{\leq Q_2}\|_2\|b_p\|_\infty\|b_q\|_2\\
\lesssim &\sum_{q> Q_2}\lambda_q^{2s+1}\|b_q\|_2\|b_q\|_\infty\sum_{p\leq Q_2}\lambda_{p}\|b_{p}\|_2\\
\lesssim &c_0m \sum_{q> Q_2}\lambda_q^{2s+1-\delta}\|b_q\|_2\Lambda_2^\delta\sum_{p\leq Q_2}\lambda_{p}\|b_{p}\|_2\\
\lesssim &c_0m \sum_{q> Q_2}\lambda_q^{s+1}\|b_q\|_2\sum_{p\leq Q_2}\lambda_{p}^{s+1}\|b_{p}\|_2\lambda_q^{s-\delta}\Lambda_2^\delta\lambda_p^{-s}.
\end{split}
\end{equation}
It follows from Young's and Jensen's inequality that
\begin{equation}\notag
\begin{split}
A_5\lesssim &c_0m \sum_{q> Q_2}\lambda_q^{2s+2}\|b_q\|_2^2+c_0m \sum_{q> Q_2}\left(\sum_{p\leq Q_2}\lambda_{p}^{s+1}\|b_{p}\|_2\lambda_q^{s-\delta}\Lambda_2^\delta\lambda_p^{-s}\right)^2\\
\lesssim &c_0m \sum_{q\geq -1}\lambda_q^{2s+2}\|b_q\|_2^2
\end{split}
\end{equation}
provided $\delta>s$.

By Lemma \ref{le-Hall1}, the definition of $\Lambda_2(t)$ \eqref{lb}, we have
\begin{equation}\notag
\begin{split}
A_{6}\lesssim & \sum_{p>Q_2+2}\sum_{|p-q|\leq 2}\lambda_q^{2s+1}\|\nabla b_{(Q_2,p-2]}\|_\infty\|b_p\|_2\|b_q\|_2\\
\lesssim & \sum_{p>Q_2+2}\|b_p\|_2\sum_{|p-q|\leq 2}\lambda_q^{2s+1}\|b_q\|_2\sum_{Q_2<p'\leq p-2}\|\nabla b_{p'}\|_\infty\\
\lesssim &c_0m \sum_{p>Q_2+2}\|b_p\|_2\sum_{|p-q|\leq 2}\lambda_q^{2s+1}\|b_q\|_2\sum_{Q_2<p'\leq p-2}\lambda_{p'}\\
\lesssim &c_0m \sum_{p>Q_2}\lambda_p^{2s+1}\|b_p\|_2^2\sum_{Q_2<p'\leq p-2}\lambda_{p'}\\
\lesssim &c_0m \sum_{q> Q_2}\lambda_q^{2s+2}\|b_q\|_2^2.
\end{split}
\end{equation} 

The term $I_{513}$ is estimated as follows,
\begin{equation}\notag
\begin{split}
|I_{513}|\leq &\sum_{q\geq-1}\sum_{|p-q|\leq 2}\ls\int_{\mathbb R^3}\left|(b_{\leq p-2}-b_{\leq q-2})\times(\nabla\times (b_p)_q)\cdot\nabla\times b_q\right|\, dx\\
\lesssim &\sum_{q>Q_2}\sum_{|p-q|\leq 2}\ls\|\nabla b_q\|_\infty\|b_{\leq p-2}-b_{\leq q-2}\|_2\|\nabla b_p\|_2\\
&+\sum_{-1\leq q\leq Q_2}\sum_{|p-q|\leq 2}\ls\|\nabla b_q\|_\infty\|b_{\leq p-2}-b_{\leq q-2}\|_2\|\nabla b_p\|_2\\
\lesssim &c_0m\sum_{q>Q_2}\sum_{|p-q|\leq 2}\lambda_q^{2s+1}\|b_{\leq p-2}-b_{\leq q-2}\|_2\|\nabla b_p\|_2\\
&+f(t)\sum_{-1\leq q\leq Q_2}\sum_{|p-q|\leq 2}\lambda_q^{2s-1}\|b_{\leq p-2}-b_{\leq q-2}\|_2\|\nabla b_p\|_2\\
\lesssim &c_0m\sum_{q\geq -1}\lambda_q^{2s+2}\| b_q\|_2^2
+f(t)\sum_{q\geq -1}\lambda_q^{2s}\| b_q\|_2^2.
\end{split}
\end{equation}
By applying Lemma \ref{le-Hall2} instead of \ref{le-Hall1}, the term $I_{52}$ can be estimated similarly as for $I_{51}$. Hence  
\[|I_{52}|\lesssim c_0m \sum_{q\geq -1}\lambda_q^{2s+2}\|b_q\|_2^2+Q_2f(t)\sum_{ q\geq -1}\lambda_q^{2s}\|b_q\|_2^2.\]
To estimate $I_{53}$, we proceed as  
\begin{equation}\notag
\begin{split}
|I_{53}|\leq &\sum_{q\geq -1}\sum_{p\geq q-2}\lambda_q^{2s}\int_{\mathbb R^3}|\Delta_q(b_p\times \nabla\times\tilde b_p)\cdot\nabla\times b_q|\, dx\\
\lesssim &\sum_{q> Q_2}\lambda_q^{2s+1}\|b_q\|_\infty\sum_{p\geq q-3}\lambda_p\|b_p\|_2^2
+\sum_{-1\leq q\leq Q_2}\lambda_q^{2s+1}\|b_q\|_\infty\sum_{p\geq q-3}\lambda_p\|b_p\|_2^2\\
\lesssim &c_0m\sum_{q> Q_2}\lambda_q^{2s+1}\sum_{p\geq q-3}\lambda_p\|b_p\|_2^2+f(t)\sum_{-1\leq q\leq Q_2}\lambda_q^{2s-1}\sum_{p\geq q-3}\lambda_p\|b_p\|_2^2\\
\lesssim &c_0m\sum_{p\geq Q_2-3}\lambda_p^{2s+2}\|b_p\|_2^2\sum_{Q_2< q\leq p+3}\lambda_{q-p}^{2s+1}
+f(t)\sum_{p\geq -1}\lambda_p^{2s}\|b_p\|_2^2\sum_{q\leq p+3}\lambda_{q-p}^{2s-1}\\
\lesssim &c_0m\sum_{q\geq Q_2-3}\lambda_q^{2s+2}\|b_q\|_2^2+f(t)\sum_{q\geq -1}\lambda_q^{2s}\|b_q\|_2^2,
\end{split}
\end{equation}
provided $s>\frac12$.
Therefore, for $s>\frac12$, we have 
\begin{equation}\label{est-i5}
|I_5|\lesssim c_0m\sum_{q\geq -1}\lambda_q^{2s+2}\|b_q\|_2^2+Q_2f(t)\sum_{q\geq -1}\lambda_q^{2s}\|b_q\|_2^2.
\end{equation}

Recall that $m=\min\{\nu,\mu\}$. Combining (\ref{est-i1})--(\ref{est-i5}), there exist absolute constants $C_1$ and $C_2$ such that
\begin{equation}\label{est-is}
\begin{split}
&|I_1|+|I_2|+|I_3|+|I_4|+|I_5|\\
\leq &C_1c_0\min\{\nu,\mu\}\left(\| u\|_{\dot H^{s+1}}^2+\| b\|_{\dot H^{s+1}}^2\right)\\
&+C_2 \max\{Q_1, Q_2\} f(t)
\left(\| u\|_{\dot H^{s}}^2+\|b\|_{\dot H^{s}}^2\right).
\end{split}
\end{equation}
Take $c_0=\frac1{C_1}$. It then follows from (\ref{ineq-ubq}) and (\ref{est-is}) that
\begin{equation}\label{energy2}
\frac{d}{dt}\left(\| u\|_{\dot H^{s}}^2+\|b\|_{\dot H^{s}}^2\right)
\leq C_2\max\{Q_1, Q_2\} \Lambda_2(t)f(t)\left(\| u\|_{\dot H^{s}}^2+\|b\|_{\dot H^{s}}^2\right).
\end{equation}
Next we show that the factor $\max\{Q_1, Q_2\}$ does not cause a problem. Indeed, we infer from (\ref{Q}) and Bernstein's inequality that
\[\Lambda_1\leq (c_0m)^{-1}\|u_{Q_1}\|_{\infty}\leq (c_0m)^{-1}\Lambda_1^{\frac32}\|u_{Q_1}\|_2=(c_0m)^{-1}\Lambda_1^{\frac32-s}\lambda_{Q_1}^s\|u_{Q_1}\|_2.\]
Thus, it indicates
\[\Lambda_1^{s-\frac12}\leq (c_0m)^{-1}\|u\|_{\dot H^s}.\]
Since $s>\frac12$, it follows then 
\[Q_1=\log\Lambda_1\leq C(\nu,\mu, s)\left(1+\log \|u\|_{\dot H^s}\right).\]
Similarly, we can deduce that for $s>\frac32$, 
\[\Lambda_2^{s-\frac32}\leq (c_0m)^{-1}\|b\|_{\dot H^s}, \ \ Q_2=\log\Lambda_2\leq C(\nu,\mu, s)\left(1+\log \|b\|_{\dot H^s}\right).\]
Therefore, from (\ref{energy2}), we obtain
\[\frac{d}{dt}\left(\| u\|_{\dot H^{s}}^2+\|b\|_{\dot H^{s}}^2\right)
\leq C(\nu,\mu,s)f(t)\left(1+\log (\|u\|_{\dot H^s}+ \|b\|_{\dot H^s})\right)\left(\| u\|_{\dot H^{s}}^2+\|b\|_{\dot H^{s}}^2\right).\]
We conclude that, by Gr\"onwall's inequality, $\|u\|_{ \dot H^s}^2+\|b\|_{\dot H^s}^2$ is  bounded on $[0,T)$ provided $f\in L^1(0,T)$. Notice that the statement holds for any  $s>3/2$.

\subsection{Comparison}

In the following lemmas, we show that the regularity condition $f\in L^1(0,T)$ is weaker than the Prodi-Serrin type criteria and an improvement criterion in the limit case of the Prodi-Serrin type  in \cite{CL}. 

\begin{Lemma}\label{compare}
Let $(u(t),b(t))$ be a weak solution to (\ref{HMHD}) on $[0,T]$. If $u\in L^q((0,T);  B^{-\frac3p}_{\infty,\infty})$ and $\nabla b\in L^\gamma((0,T);  B^{-\frac3\beta}_{\infty,\infty})$ with $\frac2q+\frac3p=1$, $\frac2\gamma+\frac3\beta=1$, and $q, \gamma\geq 2$ and $p, \beta> 3$, then $f\in L^1(0,T)$.
\end{Lemma}
\pf 
Let $f_1(t)= \|u_{\leq Q_1(t)}\|_{ B^1_{\infty,\infty}}$ and $f_2=\Lambda_2(t)\|b_{\leq Q_2(t)}\|_{ B^1_{\infty,\infty}}$. Following the lines of proof of Lemma 4.2 in \cite{CSreg}, we can show $f_1, f_2\in L^1(0,T)$. Hence $f=f_1+f_2\in L^1(0,T)$.
\cbdu

In \cite{CL}, the authors proved that if a regular solution $(u,b)$ on $[0,T)$ satisfies \eqref{cl-criterion1}--\eqref{cl-criterion2} or \eqref{cl-bmo}
then the regular solution can be extended beyond the time $T$.
Notice that the following embedding holds
\[L^\gamma((0,T); L^\beta)\subset L^\gamma((0,T);  B^0_{\beta,\infty})\subset L^\gamma((0,T); B^{-\frac3\beta}_{\infty,\infty}).\]
Thus, if $(u,b)$ satisfies \eqref{cl-criterion1}--\eqref{cl-criterion2}, then $f\in L^1(0,T)$. On the other hand, due to the embedding $BMO\subset B^0_{\infty,\infty}$ and hence $L^2(0,T; BMO)\subset L^2(0,T; B^0_{\infty,\infty})$, it follows from Lemma \ref{compare} with $(q=\gamma=2, \, p=\beta=\infty)$ that if $(u,b)$ satisfies \eqref{cl-bmo}, then $f\in L^1(0,T)$.

\bigskip

\medskip

{\bf Acknowledgement.} The author would like to thank the anonymous referee for the careful review and valuable suggestions which have helped to improve the paper a lot.


\bigskip



\begin{thebibliography}{XX}



\bibitem{ADFL}
M. Acheritogaray, P. Degond, A. Frouvelle and J-G. Liu.
\newblock {\em Kinetic formulation and global existence for the Hall-Magnetohydrodynamic system}.
\newblock Kinetic and Related Models, 4: 901--918, 2011.

\bibitem{BCD}
H. Bahouri, J. Chemin,  and R. Danchin.
\newblock {\em Fourier analysis and nonlinear partial differential equations}.
\newblock Grundlehrender Mathematischen Wissenschaften, 343. Springer, Heidelberg, 2011.


\bibitem{CDL}
D. Chae, P. Degond and J-G. Liu.
\newblock {\em Well-posedness for Hall--magnetohydrodynamics}.
\newblock arXiv:1212.3919, 2012.

\bibitem{CL}
D. Chae and J. Lee.
\newblock {\em On the blow-up criterion and small data global existence for the Hall-magneto-hydrodynamics}.
\newblock J. Differential Equations, 256: 3835--3858, 2014.

\bibitem{CS}
D. Chae,  and M. Schonbek.
\newblock {\em On the temporal decay for the Hall-magnetohydrodynamic equations}.
\newblock arXiv:1302.4601, 2013.

\bibitem{CWW}
D. Chae,  R. Wan and J. Wu.
\newblock {\em Local well-posedness for the Hall--MHD equations with fractional magnetic diffusion}.
\newblock arXiv:1404.0486v2, 2014.

\bibitem{CW}
D. Chae and J. Wolf.
\newblock {\em On partial regularity for the 3D non-stationary Hall magnetohydrodynamics equations on the plane}.
\newblock arXiv:1502.0347, 2015.

\bibitem{CMZ}
Q. Chen, C. Miao, and Z. Zhang.
\newblock {\em A new Bernstein inequality and the 2D dissipative quasi-geostrophic equation}.
\newblock Commun. Math. Phys., Vol. 271: 821--838, 2007.

\bibitem{CD}
A. Cheskidov and M. Dai.
\newblock {\em Regularity criteria for the 3D Navier-Stokes and MHD equations}.
\newblock arXiv:1507.06611, 2015.

\bibitem{CSreg}
A. Cheskidov and R. Shvydkoy.
\newblock {\em A unified approach to regularity problems for the 3D Navier-Stokes and Euler equations: the use of Kolmogorov's dissipation range}.
\newblock J. Math. Fluid Mech., Vol.16, Issue 2: 263--273, 2014.



\bibitem{CC}
A. C\'ordoba and D. C\'ordoba.
\newblock {\em A maximum principle applied to quasi-geostrophic equations}.
\newblock Comm. Math. Phys., 249(3): 511--528, 2004.

\bibitem{D}
M. Dai.
\newblock {\em Regularity criterion and energy conservation for the supercritical Quasi-Geostrophic equation}.
\newblock Journal of Mathematical Fluid Mechanics. To appear. ArXiv:1505.02293, 2015.

\bibitem{DS}
E. Dumas and F. Sueur.
\newblock {\em On the weak solutions to the Maxwell-Landau-Lifshitz equations and to the Hall-magnetohydrodynamic equations}.
\newblock Comm. Math. Phys., 330: 1179--1225, 2014.


\bibitem{For}
T. G. Forbes.
\newblock {\em Magnetic reconnection in solar flares}.
\newblock Geophys. astropphys. fluid dynamics, 62: 15--36, 1991.


\bibitem{Gr}
L. Grafakos.
\newblock {\em Modern Fourier analysis}.
\newblock Second edition. Graduate Texts in Mathematics, 250. Springer, New York, 2009.

\bibitem{K41}
A.~Kolmogoroff.
\newblock The local structure of turbulence in incompressible viscous fluid for very large {R}eynold's numbers.
\newblock {\em C. R. (Doklady) Acad. Sci. URSS (N.S.)}, 30:301--305, 1941.

\bibitem{L}
P. G. Lemari\'e-Rieusset.
\newblock {\em Recent developments in the Navier-{S}tokes problem}.
\newblock Chapman and Hall/CRC Research Notes in Mathematics, 431. Chapman  and Hall/CRC, Boca Raton, FL, 2002.

\bibitem{Li}
M. J. Lighthill.
\newblock {\em Studies on magnetohydrodynamic waves and other anisotropic wave motions}.
\newblock Philos. Trans. R. Soc. Lond., Ser. A : 397--430, 1960.






\bibitem{SC}
A. N. Simakov and L. Chacon.
\newblock {\em Quantitative, analytical model for magnetic reconnection in Hall magnetohydrodynamics}.
\newblock Phys. Rev. Lett, 101, 105003, 2008.



\end{thebibliography}
\end{document}